\documentclass[11pt,a4paper]{article}
\usepackage{amsmath}
\usepackage{amssymb}
\usepackage{caption}
\usepackage{multirow}
\usepackage{graphicx}
\usepackage{float}
\usepackage{array}
\usepackage{subcaption}
\newtheorem{theorem}{Theorem}[section]

\DeclareMathOperator{\tr}{tr}
\numberwithin{equation}{section}
\numberwithin{theorem}{section}
\numberwithin{corollary}{section}

\begin{document}

\title{Stability of equilibriums and bifurcation analysis of  two-dimensional autonomous competitive Lotka-Volterra dynamical system}

\author{Danijela Brankovi\'c and Marija Miki\'c}

\date{}
\maketitle

{\bf Abstract.} A detailed analysis of the stability of equilibriums and bifurcations of the two-dimensional autonomous  competitive Lotka-Volterra dynamical system is performed. Necessary and sufficient conditions are determined for equilibriums (without the origin) to be asymptotically stable or unstable on $\left[0, +\infty\right)^2$. Necessary and sufficient conditions are determined so that the observed dynamical system has no equilibriums in $\left(0, +\infty\right)^2 $. All results are presented in five tables and five figures. We also found that four transcritical bifurcations occur in the observed dynamical system if it is analyzed on $\mathrm{R}^2$.
\smallskip

{\emph{2020 Mathematics Subject Classification:} 34A34; 34C23; 37C75.}
\smallskip

{\emph{Keywords and phrases:} Lotka-Volterra; dynamical systems; stability; equilibriums; bifurcations.

\section{Introduction}

In this paper, we analyze the stability of equilibriums and bifurcations of a two-dimensional autonomous competitive Lotka-Volterra dynamical system
\begin{equation}\label{2d}
\begin{array}{l}
\displaystyle \dot{x_1}=x_1 \left(b_1-a_{11}x_1-a_{12}x_2\right),\\
\displaystyle \dot{x_2}=x_2 \left(b_2-a_{21}x_1-a_{22}x_2\right),		
\end{array}
\end{equation}
where $x_i=x_i\left(t\right)$ is the density of  species $i$ at a given time $t$, $b_i$ is the inherent growth rate of species $i$, $a_{ij}$ is the effect that species $j$ has on species $i$, for $i, j \in \{1,2\}$ and $i\neq j$, while for $i=j$ we have self-interacting terms $a_{ii}$, $i \in \{1,2\}$. In the absence of competition, the observed dynamical system (\ref{2d}) consists of two logistic equations (one for each species), with carrying capacity $b_i/a_{ii}$, for $i \in \{1,2\}$. Since the dynamical system (\ref{2d}) is competitive, we assume that $a_{ij}>0$ and $b_i>0$, for every $i, j \in \{1,2\}$. The reason for this is the assumed harmfulness of interspecies interactions, which directly implies $a_{ij}>0$, for every $i, j \in \{1,2\}$ and $i \neq j$. Moreover, the self-regulation of each species implies $a_{ii}>0$, for every $i \in \{1, 2\}$. Moreover $b_i>0$, because it is also assumed that the inherent growth rate of species $i$ is positive in the absence of competition, unless $x_i$ has the value $b_i/a_{ii}$ (its carrying capacity), for $i \in \{1, 2\}$. 

\noindent Competitive dynamical systems were introduced by Lotka and Volterra (see \cite{f8, f10}). The stability of equilibriums of the autonomous Lotka-Volterra  competitive model is analyzed in various literature (see, e.g. \cite{lvspanish, murray, hofbauer}). In \cite{zeeman1995} relationship between the coefficients $a_{ij}$ and $b_i$ can be found  so that the system (\ref{2d}), which has no equilibrium in $\left(0,+\infty\right)^2$, has an asymptotically stable equilibrium on one axis, while an unstable equilibrium lies on the other axes. These conditions are also generalized for the $n$-dimensional case of the system (\ref{2d}) in \cite{zeeman1995}, and \cite{ahmadlazer1996}. Other papers dealing with suitable conditions for the two-dimensional and $n$-dimensional  non-autonomous Lotka-Volterra competition model are \cite{ahmad1993nonautonomous} and \cite{ahmadlazer1998}, respectively. 

\noindent The importance of this work lies in the fact that all possible cases concerning the stability of equilibriums of the autonomous Lotka-Volterra competitive model are discussed here in detail, which to our knowledge has not been done anywhere else in the literature in this way. The new results of this paper concern the improvement of several results from \cite{zeeman1995}, specifically the determination of necessary and sufficient conditions for equilibriums (without the origin) to be asymptotically stable or unstable in $\left[0, +\infty\right)^2$, as well as the determination of necessary and sufficient conditions for the observed dynamical system to have no equilibriums in $\left(0, +\infty\right)^2$. Furthermore, we have concluded that four transcritical bifurcations occur in the observed dynamical system if it is analyzed on $\mathrm{R}^2$, and we have determined the conditions under which these bifurcations occur. To our knowledge, these are also new results of this paper.

\noindent In Section~\ref{sectionstability} we analyze in detail the stability of equilibriums of the system (\ref{2d}). We assume that $x_i \geq 0$ and that the main determinant or the minor determinants of the system (\ref{2d}) can be zero. This leads to some non-hyperbolic equilibriums and bifurcations. We note that in \cite{florian} only the stability of equilibriums of the system (\ref{2d}) for $x_i>0$ and for all non-zero determinants of this system was studied. This analysis can also be found in our paper in  Section~\ref{nzero12}, but with a new notation that is more suitable for generalization to higher dimensions. In Section~\ref{theorems} we have improved several two-dimensional versions of the theorems from the paper \cite{zeeman1995}, in particular with respect to sufficient and necessary conditions for equilibriums of the system (\ref{2d}) lying on the axes to be asymptotically stable or  an unstable equilibriums. In Section~\ref{sectionbifurcation} we present a bifurcation analysis of the dynamical system (\ref{2d}).

\noindent First we introduce the following notation. For the  determinant of the matrix of the system (\ref{2d}) and the corresponding minors, we introduce the following notation 
\begin{equation}{\label{not}}
\begin{array}{l}
\displaystyle d_{12}= \left|\begin{matrix}
a_{11} & a_{12}\\
a_{21} & a_{22}
\end{matrix}\right|,\,\,
							
\displaystyle d_{122}= -\left|\begin{matrix}
b_1 & a_{12}\\
b_2 & a_{22}
\end{matrix}\right|,\,\,\mbox{and}\,\,\,
												
\displaystyle  d_{112}= \left|\begin{matrix}
a_{11} & b_1\\
a_{21} & b_2
\end{matrix}\right|.
\end{array}							
\end{equation}

\noindent From (\ref{not}) we can conclude that the number of different indices of the determinant indicates the order of this determinant as well as the elements on its main diagonal. The indices of the determinant $d_{12}$ of the matrix of the system (\ref{2d}) show that it is the second-order determinant whose elements on the main diagonal are $a_{11}$ and $a_{22}$. Furthermore, the indices of the minors $d_{ijk}$, where $i, j, k \in \{1, 2\}$ denote the second-order determinant whose elements on the main diagonal are $a_{ij}$ (position $\left(1,1\right)$) and $b_k$ (position $\left(2,2\right)$). This notation is very useful for calculating the coordinates of the equilibriums of the system and their eigenvalues, as well as for generalizations to higher dimensions. 

\noindent The following notation applies to Tables 1-5: U denotes an unstable equilibrium, AS represents an asymptotically stable equilibrium, SS stands for a semi-stable equilibrium and NI stands for a non-isolated equilibrium. 

\noindent In all figures in our paper, we have colored the points representing unstable equilibriums in yellow, asymptotically stable equilibriums are colored red, semi-stable equilibriums are colored orange and non-isolated equilibriums are colored pink.

\section{Stability analysis of the equilibriums}\label{sectionstability}

We consider two cases, depending on whether $d_{12}$ is zero or not.

\subsection{Case \boldmath{$d_{12} \neq 0$}}\label{nzero}

The number of equilibriums and their stability depend on the signs of $d_{12}$, $d_{112}$ and $d_{122}$. We analyze the following cases.

\subsubsection{Case \boldmath{$d_{112} \neq 0$} and \boldmath{$d_{122} \neq 0$}}\label{nzero12}

If $\mathrm{sgn}\left(d_{12}\right)=\mathrm{sgn}\left(d_{112}\right)=-\mathrm{sgn}\left(d_{122}\right)$, we have four different equilibriums of system (\ref{2d}), labelled $E_0\left(0,0\right)$, $E_1\left(b_1/a_{11},0\right)$, $E_2\left(0,b_2/a_{22}\right)$ and $E_{12}\left(-d_{122}/d_{12},d_{112}/d_{12}\right)$. Otherwise we would have three equilibriums of system (\ref{2d}), $E_0$, $E_1$ and $E_2$, since equilibrium $E_{12}$ would not lie in the first quadrant. We note that the conditions $d_{112} \neq 0$ and $d_{122} \neq 0$ ensure that the eigenvalues of the Jacobian matrix in our equilibriums $E_i$, $i \in \{1,2,12\}$ have non-zero real parts, as we will see in the following text. Consequently, our equilibriums are hyperbolic and therefore we can discuss their stability in the framework of the corresponding linearization of the dynamical system (\ref{2d}), because the Hartman--Grobman theorem provides us topological equivalence between a nonlinear dynamical system and its linearization in the neighbourhood of the hyperbolic equilibrium. Now we investigate the stability of all equilibriums. 

\noindent{\bf{The equilibrium \boldmath{$E_0$}.}} The eigenvalues of the Jacobian matrix at $E_0$ are $\lambda_i=b_i$, for $i \in \{1,2\}$. Since $b_i>0$ for every $i \in \{1,2\}$, we conclude that $E_0$ is an unstable node (unstable equilibrium). 

\noindent{\bf{The equilibrium \boldmath{$E_1$}.}} The eigenvalues of the Jacobian matrix at $E_1$ are $\lambda_1=-b_1$, $\lambda_1<0$ and $\lambda_2=d_{112}/a_{11}$. If $d_{112}<0$, then $E_1$ is a stable node (asymptotically stable equilibrium). Otherwise, i.e.\ if $d_{112}>0$, $E_1$ is therefore a saddle (unstable equilibrium).

\noindent{\bf{The equilibrium \boldmath{$E_2$}.}} We note that the eigenvalues of the Jacobian matrix for $E_2$ are $\lambda_1=-d_{122}/a_{22}$ and $\lambda_2=-b_2$, $\lambda_2<0$. We conclude that $E_2$ is a stable node (asymptotically stable equilibrium) if $d_{122}>0$. If $d_{122}<0$, then $E_2$ is a saddle point (unstable equilibrium).

\noindent{\bf{The equilibrium \boldmath{$E_{12}$}.}} The situation is a little complicated in this case. As we have already mentioned, the coordinates of $E_{12} \left(-d_{122}/d_{12}, d_{112}/d_{12} \right)$ must be positive in order to $E_{12}$ lie in the first quadrant, i.e.\ it must be $\mathrm{sgn}\left(d_{12}\right)=\mathrm{sgn}\left(d_{112}\right)=-\mathrm{sgn}\left(d_{122}\right)$. If we denote the coordinates of $E_{12}$ by $x_1:=-d_{122}/d_{12}$ and $x_2:=d_{112}/d_{12}$, the Jacobian determinant at $E_{12}$ is
	\begin{equation}\label{jac12}
		\begin{array}{ll}
		\smallskip
		\displaystyle \det \left(J\left(E_{12}\right)\right)&=\left|\begin{matrix}
											b_1-2a_{11}x_1-a_{12}x_2 & -a_{12}x_1\\
											-a_{21}x_2 & b_2-a_{21}x_1-2a_{22}x_2 
											\end{matrix}\right|\\
											
											\smallskip
											&=
											\left(-1\right)^2 x_1x_2\left|\begin{matrix}
																	a_{11} & a_{12}\\
																	a_{21} & a_{22} 
																														\end{matrix}\right|=x_1 x_2 d_{12},
																					\end{array}	
	\end{equation}	
since $b_1-2a_{11}x_1-a_{12}x_2=-a_{11}x_1$ and $ b_2-a_{21}x_1-2a_{22}x_2=-a_{22}x_2$.
Furthermore
	\begin{equation}\label{pte}
		\begin{array}{l}
		\displaystyle \lambda_1 \lambda_2 = \det\left( J\left(E_{12}\right)\right)=x_1 x_2 d_{12},\\
		
		\smallskip
		\displaystyle \lambda_1+\lambda_2 = \tr \left( J\left(E_{12}\right)\right)=-\left(a_{11}x_1+a_{22}x_2\right),		
		\end{array}	
	\end{equation}
from which we obtain that the characteristic polynomial $P_{12}\left(\lambda\right)$ of the matrix $J\left(E_{12}\right)$ is $P_{12}\left(\lambda\right)=\lambda^2+\left(x_1 a_{11}+x_2 a_{22}\right)\lambda + x_1 x_2 d_{12}$. From (\ref{pte}) we conclude that $E_{12}$ is a stable node (asymptotically stable equilibrium) if $d_{12}>0$ and $d_{112}>0$ and $d_{122}<0$. If $d_{12}<0$ and $d_{112}<0$ and $d_{122}>0$, $E_{12}$ is a saddle (unstable equilibrium).

\noindent To summarize, $E_0$ is an unstable node (unstable equilibrium), while the stability of equilibrium $E_1$ depends solely on the sign of $d_{112}$, the stability of equilibrium $E_2$ depends solely on the sign of $d_{122}$, while the existence and stability of equilibrium $E_{12}$ depend on the signs of $d_{12}$, $d_{112}$ and $d_{122}$.  The results are shown in Table~\ref{Table1}. 

	\begin{table}[h!]
		\centering
			\begin{tabular}{||>{\centering\arraybackslash}m{1.1cm}||>{\centering\arraybackslash}m{0.4cm}| >{\centering\arraybackslash}m{2.1cm}|>{\centering\arraybackslash}m{2.1cm}|>{\centering\arraybackslash}m{2.1cm}|>{\centering\arraybackslash}m{2.1cm}||}
		\hline
	$d_{12}$ & eq. &  $d_{112}>0$ \newline $d_{122}>0$  & $d_{112}>0$ \newline $d_{122}<0$ & $d_{112}<0$ \newline $d_{122}>0$ & $d_{112}<0$ \newline $d_{122}<0$\\
		\hline\hline
		\multirow{4}{1.3cm}{$d_{12}>0$} & $E_0$ & U (u.\ node) & U (u.\ node) & \multirow{4}{1.9cm}{\centering not possible case} & U (u.\ node)\\ \cline{2-4}\cline{6-6}
		
		& $E_1$ & U (saddle) & U (saddle) &  & AS (s.\ node) \\ \cline{2-4}\cline{6-6}
		
		& $E_2$ & AS (s.\ node) & U (saddle) &  & U (saddle) \\ \cline{2-4}\cline{6-6}
		& $E_{12}$ &  not exist & AS (s.\ node) &  & not exist \\ \cline{2-4}
		
		\hline\hline
	
		\multirow{4}{1.3cm}{$d_{12}<0$} & $E_0$ & U (u.\ node) & \multirow{4}{1.9cm}{\centering not possible case} & U (u.\ node) & U (u.\ node) \\ \cline{2-3}\cline{5-6}
		
		& $E_1$ & U (saddle) &  & AS (s.\ node) & AS (s.\ node) \\ \cline{2-3}\cline{5-6} 
		
		& $E_2$ & AS (s.\ node) &  & AS (s.\ node) & U (saddle) \\ \cline{2-3}\cline{5-6}
		& $E_{12}$ &  not exist & &  U (saddle) & not exist \\ \cline{2-3}
		\hline
			\end{tabular}
		\caption{Stability analysis of the equilibriums of the system (\ref{2d}).}
		\label{Table1}
	\end{table}	
\noindent The case from Table~\ref{Table1}, in which $d_{12}>0$, $d_{112}<0$, $d_{122}>0$ is not possible. If this were possible, then we would have 
	\begin{equation}\label{impossible1}
		a_{11}a_{22}>a_{12}a_{21},\, a_{11}b_2<a_{21}b_1,\,
		a_{12}b_2>a_{22}b_1.
	\end{equation}
If we multiply the first inequality from (\ref{impossible1}) by $b_2/a_{22}$ and use the third inequality from (\ref{impossible1}), we get $ a_{11}b_2>a_{12}b_2 a_{21}/a_{22}>a_{21}b_1$, which contradicts the second inequality from (\ref{impossible1}). The case that $d_{12}<0$, $d_{112}>0$, $d_{122}<0$ is also not possible. 

\noindent We conclude from Table~\ref{Table1} that cases with $d_{12}>0$, $d_{112}>0$, $d_{122}>0$ and $d_{12}<0$, $d_{112}>0$, $d_{122}>0$, have the same qualitative dynamical properties, as well as the cases $d_{12}>0$, $d_{112}<0$, $d_{122}<0$ and $d_{12}<0$, $d_{112}<0$, $d_{122}<0$. Therefore, in this section we have four qualitatively different phase portraits of the system (\ref{2d}). They are shown in Figures~\ref{f1} and \ref{f2}.
\begin{figure}[H]
		\begin{subfigure}{0.44\textwidth}
		 	\includegraphics[width=\linewidth]{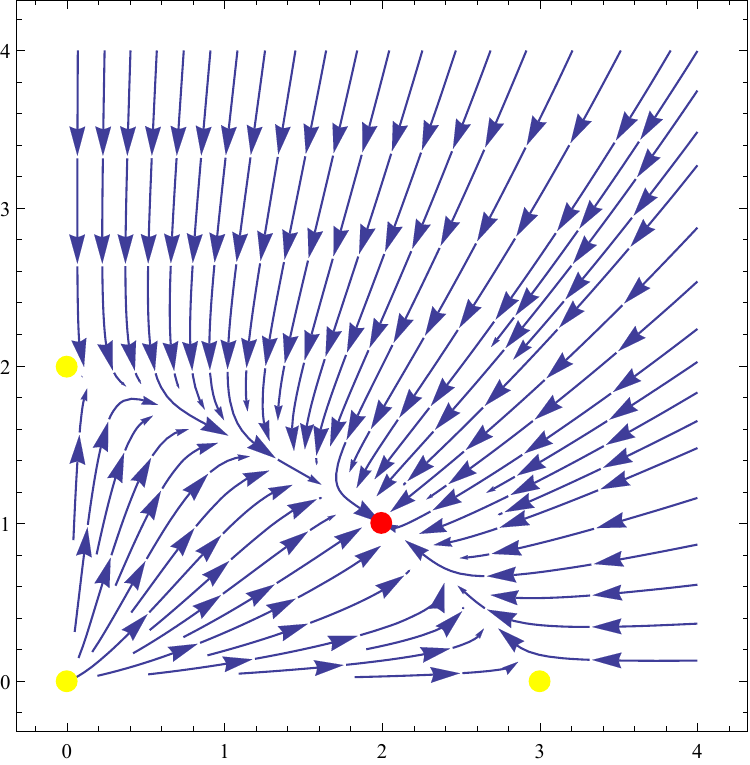}
		 	\caption{Case $d_{12}>0$, $d_{112}>0$, $d_{122}<0$. $\begin{array}{ll} \dot{x_1}&=x_1 \left(3-x_1-x_2\right), \\ \dot{x_2}&=x_2 \left(4-x_1-2x_2\right). \end{array}$}
			 \label{f1a}
		\end{subfigure}
		\hfill
		\begin{subfigure}{0.44\textwidth}
		 	\includegraphics[width=\linewidth]{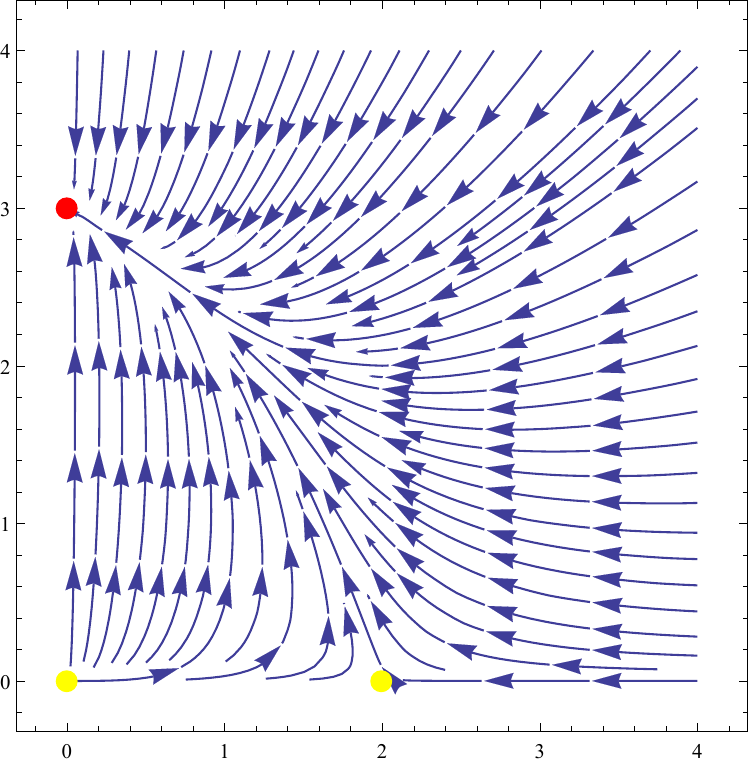}
		 	\caption{Case $d_{12}>0$, $d_{112}>0$, $d_{122}>0$. $\begin{array}{ll} \dot{x_1}&=x_1 \left(2-x_1-x_2\right), \\ \dot{x_2}&=x_2 \left(6-x_1-2x_2\right). \end{array}$}
			 \label{f1b}
		\end{subfigure}
		\hfill
		 \caption{Phase portraits for the cases $d_{12}>0$, $d_{112}>0$, $d_{122}<0$ and $d_{12}>0$, $d_{112}>0$, $d_{122}>0$, respectively.}
		 \label{f1}
	\end{figure} 
\noindent In Figures~\ref{f1a} and \ref{f1b}, the equilibrium $E_0$ at the origin is an unstable node, while the equilibrium $E_1$ is a saddle. In Figures~\ref{f1a} and \ref{f2a}, equilibrium $E_2$ is a saddle, while in Figures~\ref{f1b} and \ref{f2b} it is a stable node. Equilibrium $E_{12}$ in Figure~\ref{f1a} is a stable node, while in Figure~\ref{f2b} is a saddle.	
 \begin{figure}[H]
		\begin{subfigure}{0.44\textwidth}
		 	\includegraphics[width=\linewidth]{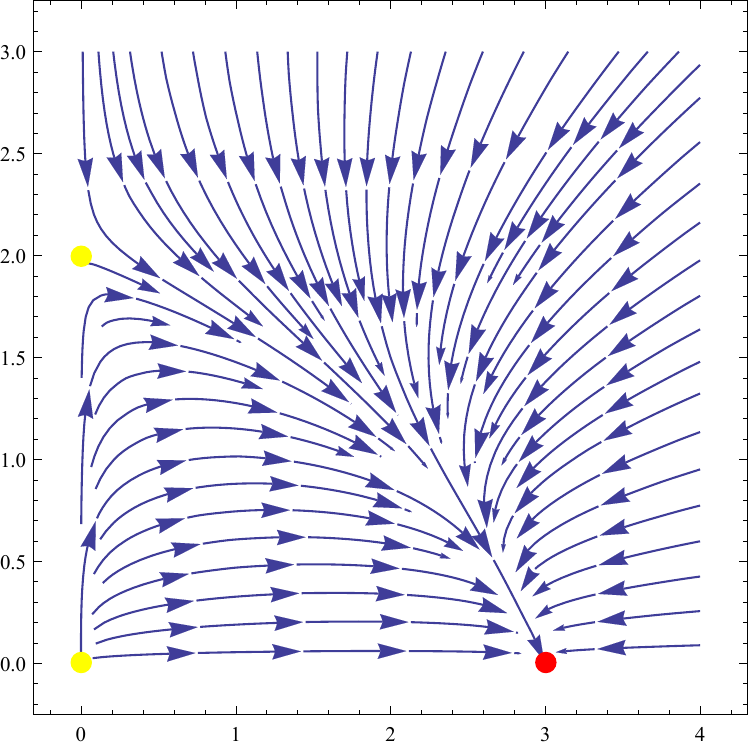}
		 	\caption{Case $d_{12}>0$, $d_{112}<0$, $d_{122}<0$. $\begin{array}{ll} \dot{x_1}&=x_1 \left(6-2x_1-x_2\right), \\ \dot{x_2}&=x_2 \left(2-x_1-x_2\right). \end{array}$}
			 \label{f2a}
		\end{subfigure}
		\hfill
		\begin{subfigure}{0.44\textwidth}
		 	\includegraphics[width=\linewidth, height=1.00\linewidth]{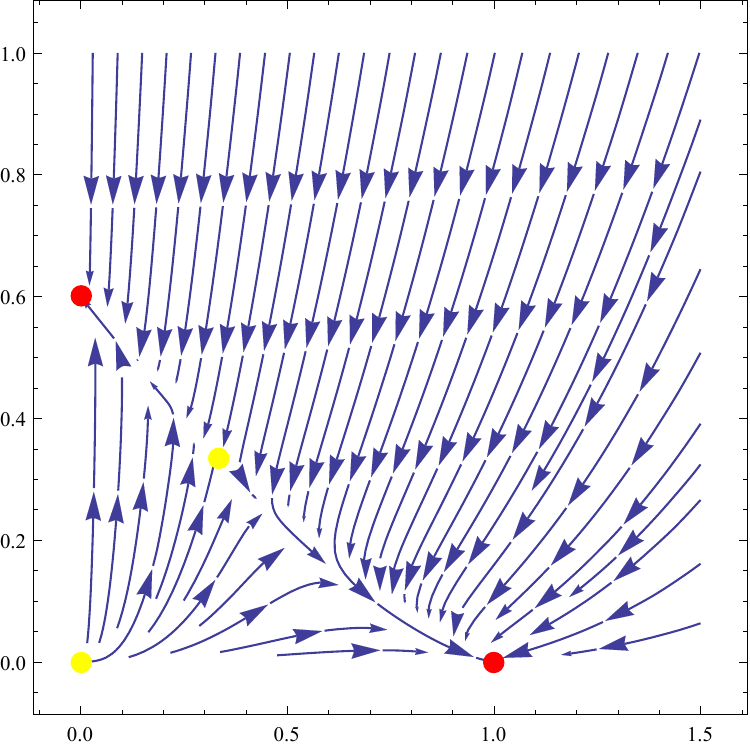}
		 	\caption{Case $d_{12}<0$, $d_{112}<0$, $d_{122}>0$. $\begin{array}{ll} \dot{x_1}&=x_1 \left(1-x_1-2x_2\right), \\ \dot{x_2}&=x_2 \left(3-4x_1-5x_2\right). \end{array}$}
			 \label{f2b}
		\end{subfigure}
		\hfill
		 \caption{Phase portraits for the cases $d_{12}>0$, $d_{112}<0$, $d_{122}<0$ and $d_{12}<0$, $d_{112}<0$, $d_{122}>0$, respectively.}
		 \label{f2}
	\end{figure} 

\subsubsection{Case \boldmath{$d_{112} \neq 0$} and \boldmath{$d_{122} = 0$} }\label{nzero2}

Since $d_{122} = 0$, we directly conclude that the first coordinate of $E_{12}$ is zero and that $a_{12}=a_{22}b_1/b_2$. Substituting $a_{12}=a_{22}b_1/b_2$ into $d_{12}$, we conclude that $d_{12}=a_{22}d_{112}/b_2$. Furthermore, the second coordinate of $E_{12}$ is 
$d_{112}/d_{12}=b_2/a_{22}$, $b_2/a_{22}>0$, from which it follows that the equilibrium $E_{12}$ coincides with $E_2$ and that $\mathrm{sgn}\left(d_{112}\right)=\mathrm{sgn}\left(d_{12}\right)$. The equilibriums $E_0$ and $E_1$ are hyperbolic and their nature is the same as in Section~\ref{nzero12} (see Table~\ref{Table1}), considering only the cases where $\mathrm{sgn}\left(d_{112}\right)=\mathrm{sgn}\left(d_{12}\right)$. Since the equilibrium $E_2$ is not hyperbolic ($\lambda_1=-d_{122}/a_{22}=0$), we cannot determine its nature using the Hartman-Grobman theorem. Instead, we will determine its stability by analyzing the phase trajectories around this equilibrium. For this purpose, we first determine the nullclines of the system (\ref{2d}). 

\noindent We obtain that $x_1$-nullclines are the lines $x_1=0$ and $x_2=\left(b_1-a_{11}x_1\right)/a_{12}$, while $x_2$-nullclines are $x_2=0$ and $x_2=\left(b_2-a_{21}x_1\right)/a_{22}$. Note that $x_1$-nullcline $x_1=0$ and $x_2$-nullcline $x_2=0$ are invariant. The vector field through $x_1$-nullcline $x_1=0$ is vertical. If we also substitute $x_1=0$ into the second equation in (\ref{2d}), we obtain $\dot{x_2}=x_2 \left(b_2-a_{22}x_2\right)$, from which we conclude that the direction of the vector field through the nullcline $x_1=0$ is vertically upwards ($\uparrow$) for $0<x_2<b_2/a_{22}$, while for $x_2>b_2/a_{22}$ the vector field is vertically downwards ($\downarrow$).

\noindent The vector field is horizontal through $x_2$-nullcline $x_2=0$. If we insert $x_2=0$ into the first equation of the system (\ref{2d}), we obtain $\dot{x_1}=x_1 \left( b_1-a_{11}x_1\right)$. Consequently, the vector field is directed horizontally to the right ($\rightarrow$) by the nullcline $x_2=0$ for $0<x_1<b_1/a_{11}$, while for $x_1>b_1/a_{11}$ the vector field is directed horizontally to the left ($\leftarrow$).

\noindent Using that $d_{122}=0$ and $x_1>0$, we can easily see that the $x_1$-nullcline $x_2=\left(b_1-a_{11}x_1\right)/a_{12}$ is above the $x_2$-nullcline $x_2=\left(b_2-a_{21}x_1\right)/a_{22}$ if and only if $d_{12}<0$. We determine that the intersection point of these two nullclines is $E_2\left(0, b_1/a_{12}\right)$. First, we analyze the nature of the equilibrium $E_2$ when $d_{12}<0$. 

\noindent Since we conclude that $\dot{x_1}=0$ for $x_2=\left(b_1-a_{11}x_1\right)/a_{12}$, it follows that the vector field through the nullcline $x_2=\left(b_1-a_{11}x_1\right)/a_{12}$ is vertical. If we substitute $x_2=\left(b_1-a_{11}x_1\right)/a_{12}$ into the second equation in system (\ref{2d}) (the substitution only applies to the term $x_2$ in brackets), we get $\dot{x_2}=x_1 x_2 d_{12}/a_{12}$. Therefore, the direction of the vector field through the nullcline $x_2=\left(b_1-a_{11}x_1\right)/a_{12}$ is vertically downwards ($\downarrow$).

\noindent The vector field through the nullcline $x_2=\left(b_2-a_{21}x_1\right)/a_{22}$ is horizontal. If we insert $x_2=\left(b_2-a_{21}x_1\right)/a_{22}$ into the first equation of the system (\ref{2d}), we derive $\dot{x_1}=-{x_1}^2 d_{12}/a_{22}$, from which we conclude that the direction of the vector field is directed horizontally to the right ($\rightarrow$) by 
the nullcline $x_2=\left(b_2-a_{21}x_1\right)/a_{22}$. 
	
\noindent The sketch of the direction of the vector field through the nullclines implies that the equilibrium $E_2$ is unstable, which we will prove.

\noindent Consider the subset 
$$\mathcal{G}_1=\{\left(x_1, x_2\right)\in \mathbb{R}^2 \vert \left(b_2-a_{21}x_1\right)/a_{22} <x_2< \left(b_1-a_{11}x_1\right)/a_{12},\, x_2>0\}.$$ 
We note that $\mathcal{G}_1$ is invariant, i.e.\ any trajectory that starts at some point of $\mathcal{G}_1$ remains in $\mathcal{G}_1$. Let $M\left(x_1, x_2\right)$ be an arbitrary point from the area $\mathcal{G}_1$. Therefore, $x_1>0$ and $b_1-a_{11}x_1-a_{12}x_2>0$ and from the first equation of the system (\ref{2d}) we conclude that $\dot{x_1}>0$, i.e.\ $x_1=x_1\left(t\right)$ increases over time for any point in $\mathcal{G}_1$. Consequently, there is a neighbourhood $\mathcal{O}_{E_2}$ of the equilibrium $E_2$ such that for every point $M\left(x_1, x_2\right) \in \mathcal{O}_{E_2} \cap \mathcal{G}_1$ the trajectory starting at $M\left(x_1, x_2\right)$ never returns to $\mathcal{O}_{E_2} \cap \mathcal{G}_1$. We conclude that $E_2\left(0, b_2/a_{22}\right)$ is an unstable equilibrium for $d_{12}<0$.
	
\noindent For $d_{12}>0$ we will prove that the equilibrium $E_2$ is semi-stable. More precisely, we will prove that $E_2$ is asymptotically stable in $\left(0, +\infty\right)^2$ using Lyapunov stability theorem and that it is unstable in $\left(-\infty, 0\right) \times \left(0,+\infty \right)$ by analyzing trajectories around $E_2$. 

\noindent Let suppose that the Lyapunov function has the form 
	\begin{equation}\label{lyapunov1}
		\begin{array}{l}
			\displaystyle V\left(x_1, x_2\right) = {x_1}^{a_{22}}{x_2}^{-a_{12}},\\
			
			\displaystyle \dot{V}\left(x_1, x_2\right) = -d_{12}{x_1}^{a_{22}+1}{x_2}^{-a_{12}}.	
		\end{array}
	\end{equation}	
From (\ref{lyapunov1}) we conclude that $V$ is differentiable on $\left(0, +\infty \right)^2$, $V\left(E_2\right)=0$ and $V\left(x_1, x_2\right)>0$ for every $\left(x_1, x_2\right) \in \left(0, +\infty \right)^2$. Since $d_{12}>0$, then $\dot{V}\left(x_1, x_2\right)<0$ for every $\left(x_1, x_2\right) \in \left(0, +\infty \right)^2$. Therefore, the function $V$ from (\ref{lyapunov1}) is a Lyapunov function for $E_2$ and according to the Lyapunov stability theorem, $E_2$ is an asymptotically stable equilibrium on $\left(0, +\infty \right)^2$.

\noindent Consider the subset 
$$\mathcal{G}_2=\{\left(x_1, x_2\right) \in \mathbb{R}^2 \vert \left(b_2-a_{21}x_1\right)/a_{22} <x_2< \left(b_1-a_{11}x_1\right)/a_{12}\}.$$
We note that, since here $d_{12}>0$ and $x_1<0$, the line $x_2=\left(b_2-a_{21}x_1\right)/a_{22}$ is indeed below the line $x_2=\left(b_1-a_{11}x_1\right)/a_{12}$. Moreover, we notice that $x_2>b_2/a_{22}$, where $b_2/a_{22}=d_{112}/d_{12}$ is exactly the $x_2$-coordinate of the equilibrium $E_2$. 

\noindent Similarly to the case $d_{12}<0$, we find that the direction of the vector field of (\ref{2d}) by the nullcline $x_2=\left(b_1-a_{11}x_1\right)/a_{12}$ is vertically upwards ($\uparrow$) if and only if $x_1x_2>0$. The vector field by the nullcline $x_2=\left(b_2-a_{21}x_1\right)/a_{22}$ is directed horizontally to the left ($\leftarrow$). 

\noindent The directions of the vector field through the nullclines suggest that the equilibrium $E_2$ is unstable in $\mathcal{G}_2$, which we will prove. 

\noindent Let $M\left(x_1, x_2\right)$ be an arbitrary point from the area $\mathcal{G}_2$. Therefore, $x_1<0$ and $x_2< \left(b_1-a_{11}x_1\right)/a_{12}$, so we conclude from the first equation of the system (\ref{2d}) that $\dot{x_1}<0$, i.e.\ $x_1=x_1\left(t\right)$ decreases with time for every point in $\mathcal{G}_2$. Consequently, there exists neighbourhood $\mathcal{O}_{E_2}$ of the equilibrium $E_2$ such that for every point $M\left(x_1, x_2\right) \in \mathcal{O}_{E_2} \cap \mathcal{G}_2$ a  trajectory starting at  $M\left(x_1, x_2\right)$ never returns to $\mathcal{O}_{E_2} \cap \mathcal{G}_2$. We conclude that $E_2\left(0, b_2/a_{22}\right)$ is an unstable equilibrium in $\left(-\infty, 0 \right)\times \left(0, +\infty \right) $. Overall, $E_2$ is a semi-stable equilibrium.

\noindent In this section we have two cases that are not possible. The first is for $d_{12}>0$, $d_{112}<0$, $d_{122}=0$. If we assume that it is possible, then 
\begin{equation}\label{impossible3}
a_{11}a_{22}>a_{12}a_{21}, \, a_{11}b_2<a_{21}b_1,\, a_{12}b_2=a_{22}b_1.
\end{equation} 

\noindent If we multiply the first inequality from (\ref{impossible3}) by $b_2/a_{22}$ and apply the third inequality from (\ref{impossible3}), we obtain, as before, $a_{11}b_2>a_{12}b_2 a_{21}/a_{22}=a_{21}b_1$, which contradicts the second inequality from (\ref{impossible3}). The case \newline $d_{12}<0$, $d_{112}>0$, $d_{122}=0$ is also not possible.

\begin{table}[h!]
		\centering
			\begin{tabular}{||>{\centering\arraybackslash}m{1.3cm}||>{\centering\arraybackslash}m{1.0cm}||>{\centering\arraybackslash}m{3.8cm}|>{\centering\arraybackslash}m{3.8cm}||}
		\hline
		$d_{12}$ & eq. & $d_{112}>0$ and $d_{122}=0$  & $d_{112}<0$ and $d_{122}=0$\\
		\hline\hline
		 \multirow{3}{1.3cm}{$d_{12}>0$} & $E_0$ & U (u.\ node) & \multirow{3}{3.8cm}{\centering not possible case}  \\	\cline{2-3}
		 & $E_1$ & U (saddle)  &    \\ \cline{2-3}		
		 & $E_2$ & SS  &   \\

		\hline\hline
	
		 \multirow{3}{1.3cm}{$d_{12}<0$}  & $E_0$ & \multirow{3}{3.8cm}{\centering not possible case} & U (u.\ node) \\	 \cline{2-2}	 \cline{4-4}
		 &  $E_1$ &   & AS (s.\ node) \\ \cline{2-2}	 \cline{4-4}		
		 &  $E_2$ &  & U \\ 
		
		\hline
			\end{tabular}
		\caption{Stability analysis of the equilibriums of the system (\ref{2d}).}
		\label{Table2}
	\end{table}

\noindent All results are listed in Table~\ref{Table2}. To conclude this chapter, we have two different phase portraits of the dynamical system (\ref{2d}), shown in Figure~\ref{f3}.

 \begin{figure}[H]
		\begin{subfigure}{0.44\textwidth}
		 		\includegraphics[width=\linewidth]{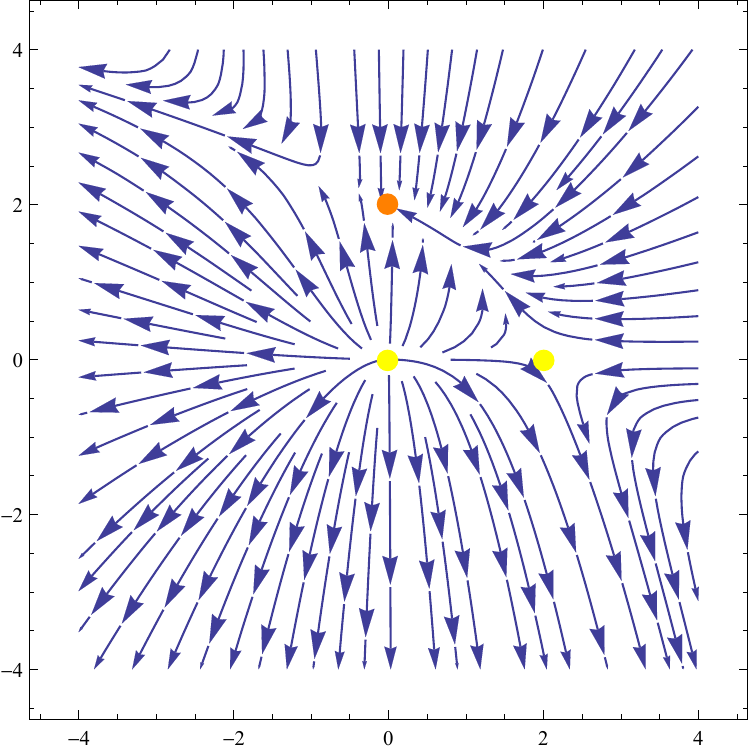}
		 	\caption{Case $d_{12}>0$, $d_{112}>0$, $d_{122}=0$. $\begin{array}{ll} \dot{x_1}&=x_1 \left(2-x_1-x_2\right), \\ \dot{x_2}&=x_2 \left(4-x_1-2x_2\right). \end{array}$}
			 \label{f3a}
		\end{subfigure}
		\hfill
		\begin{subfigure}{0.44\textwidth}
		 	\includegraphics[width=\linewidth]{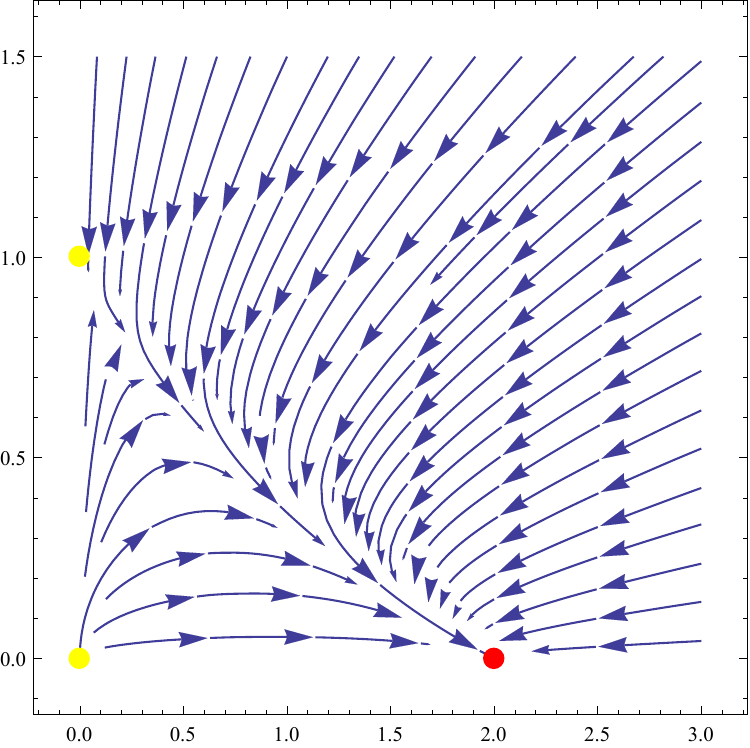}
		 	\caption{Case $d_{12}<0$, $d_{112}<0$, $d_{122}=0$. $\begin{array}{ll} \dot{x_1}&=x_1 \left(2-x_1-2x_2\right), \\ \dot{x_2}&=x_2 \left(1-x_1-x_2\right). \end{array}$}
			 \label{f3b}
		\end{subfigure}
		\hfill
		 \caption{Phase portraits for the cases $d_{12}>0$, $d_{112}>0$, $d_{122}=0$ and $d_{12}<0$, $d_{112}<0$, $d_{122}=0$, respectively.}
		 \label{f3}
	\end{figure} 
	
\noindent In both figures, the equilibrium $E_0$ at the origin is an unstable node, while in Figure~\ref{f3a} the equilibrium $E_1$ is a saddle point and the equilibrium $E_2$ is a semi-stable, while in Figure~\ref{f3b} the equilibrium $E_1$ is asymptotically stable and the equilibrium $E_2$ is an unstable equilibrium.

\subsubsection{Case \boldmath{$d_{112} = 0$ and \boldmath{$d_{122} \neq 0$} }}\label{nzero1}

Similarly as in the Section \ref{nzero2}, from $d_{112} = 0$ we conclude that the equilibrium $E_{12}$ coincides with $E_1$ and that $\mathrm{sgn}\left(d_{122}\right)=-\mathrm{sgn}\left(d_{12}\right)$.
The equilibriums $E_0$ and $E_2$ are hyperbolic and their nature is the same as in Section~\ref{nzero2} (Table~\ref{Table1}), considering only the cases when $\mathrm{sgn}\left(d_{122}\right)=-\mathrm{sgn}\left(d_{12}\right)$. The equilibrium $E_1$ is nonhyperbolic ($\lambda_2=d_{112}/a_{11}=0$), and we will determine its stability by analyzing the phase trajectories around this equilibrium. However, we will only present the final results here, since the procedure is similar to the one presented in Section~\ref{nzero2}.

\noindent If $d_{12}<0$, then equilibrium $E_1$ will be unstable (proof of this statement can be obtained using nullclines, similarly as in Section~\ref{nzero2}). If $d_{12}>0$, then the equilibrium $E_1$ is semi-stable. The equilibrium $E_1$ is asymptotically stable on $\left(0, +\infty \right)^2$. To prove this statement, the reader can use the Lyapunov function $\displaystyle V\left(x_1, x_2\right) = {x_1}^{-a_{21}}{x_2}^{a_{11}}$. The equilibrium $E_1$ is unstable in $\left(0, +\infty \right) \times \left(-\infty, 0 \right)$. To prove this statement, the reader can use nullclines, similar to Section~\ref{nzero2}. The results are shown in the Table~\ref{Table3}. 
 \begin{table}[h!]
 		\centering
 			\begin{tabular}{||>{\centering\arraybackslash}m{1.3cm}||>{\centering\arraybackslash}m{1.0cm}||>{\centering\arraybackslash}m{3.8cm}|p{3.8cm}||}
 		\hline
 		 $d_{12}$ & eq. & $d_{112}=0$ and $d_{122}>0$  & $d_{112}=0$ and $d_{122}<0$\\
 		\hline\hline
 		 \multirow{3}{1.3cm}{$d_{12}>0$} & $E_0$ & \multirow{3}{3.8cm}{\centering not possible case} & U (u.\ node)  \\	 \cline{2-2}\cline{4-4}
 		&  $E_1$ &    & SS  \\	\cline{2-2}\cline{4-4}	
 		&  $E_2$ &   & U (saddle) \\ \cline{2-2}\cline{4-4}

 		\hline\hline
	
 		 \multirow{3}{1.3cm}{$d_{12}<0$} & $E_0$ & U (u.\ node)  & \multirow{3}{3.8cm}{\centering not possible case} \\	\cline{2-3}	
 		 & $E_1$ & U  &   \\		\cline{2-3}
 		 & $E_2$ & AS  &   \\ \cline{2-3}
		
 		\hline
 			\end{tabular}
 		\caption{Stability analysis of the equilibriums of the system (\ref{2d}).}
 		\label{Table3}
  	\end{table}		
\noindent Therefore, in this section we have two different phase portraits of the dynamical system (\ref{2d}), which are shown in Figure~\ref{f4}.

 \begin{figure}[H]
 	\begin{subfigure}{0.44\textwidth}
		 	\includegraphics[width=\linewidth]{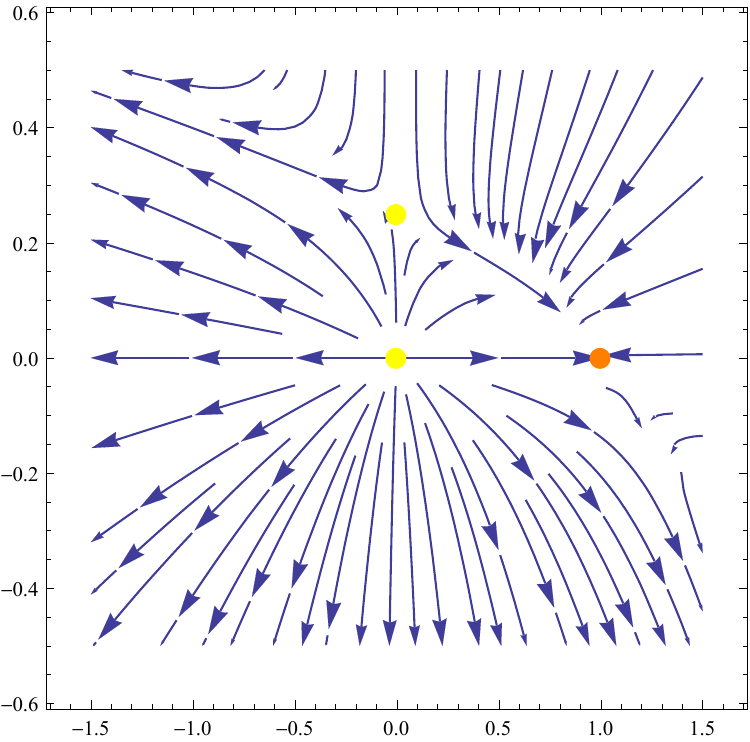}
		 	\caption{Case $d_{12}>0$, $d_{112}=0$, $d_{122}<0$. $\begin{array}{ll} \dot{x_1}&=x_1 \left(1-x_1-2x_2\right), \\ \dot{x_2}&=x_2 \left(1-x_1-4x_2\right). \end{array}$}
			 \label{f4a}
		\end{subfigure}
		\hfill
		\begin{subfigure}{0.44\textwidth}
		 		\includegraphics[width=\linewidth, height=0.98\linewidth]{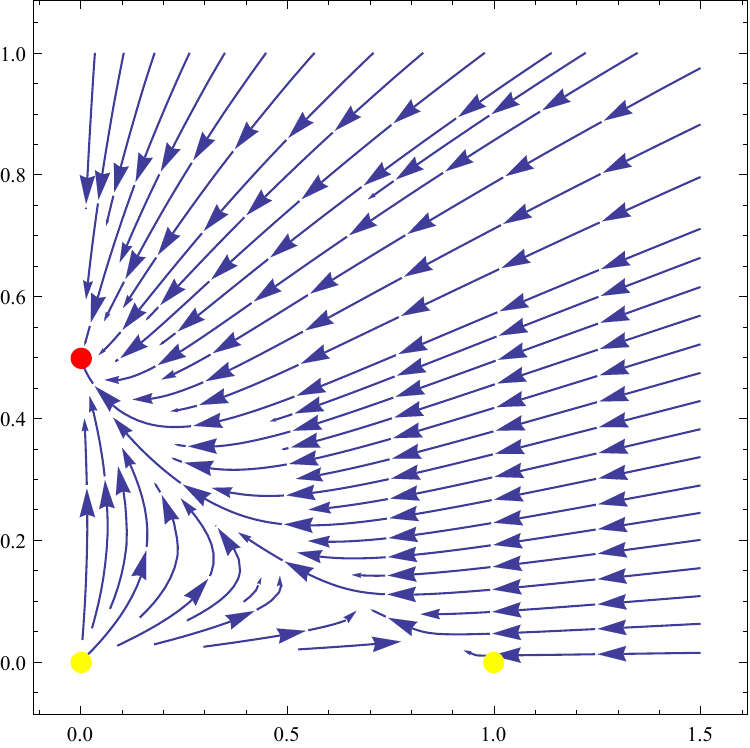}
		 	\caption{Case $d_{12}<0$, $d_{112}=0$, $d_{122}>0$. $\begin{array}{ll} \dot{x_1}&=x_1 \left(1-x_1-4x_2\right), \\ \dot{x_2}&=x_2 \left(1-x_1-2x_2\right). \end{array}$}
			 \label{f4b}
		\end{subfigure}
		\hfill
		 \caption{Phase portraits for the cases $d_{12}>0$, $d_{112}=0$, $d_{122}<0$ and $d_{12}<0$, $d_{112}=0$, $d_{122}>0$, respectively.}
		 \label{f4}
	\end{figure} 
\noindent In both figures, the equilibrium $E_0$ at the origin is an unstable node, while in Figure~\ref{f4a}, the equilibrium $E_1$ is a semi-stable and the equilibrium $E_2$ is a saddle, while in Figure~\ref{f4b} the equilibrium $E_1$ is an unstable equilibrium and the equilibrium $E_2$ is a stable node.
\subsubsection{Case \boldmath{$d_{122} = 0$} and \boldmath{$d_{112} = 0$}}\label{nzero0}

This case is not possible. This can be shown in a similar way as in Section~\ref{nzero2} (the proof differs if $d_{12}>0$ or $d_{12}<0$).  

\subsection{Case \boldmath{$d_{12} = 0$}}\label{zero}

We distinguish between two cases, when at least one minor is not equal to zero and when both minors are equal to zero.

\subsubsection{Case \boldmath{$d_{122} \neq 0$} or \boldmath{$d_{112} \neq 0$}}\label{zero1or2}

The equilibrium $E_{12}$ does not exist in this case, since the system
	\begin{equation}\label{2ds}
		\begin{array}{l}
		\smallskip
		\displaystyle a_{11}x_1+a_{12}x_2=b_1,\\
		
		\smallskip
		\displaystyle a_{21}x_1+a_{22}x_2=b_2,		
		\end{array}
	\end{equation}
has no solutions. Accordingly, the system (\ref{2d}) has three equilibriums, $E_0$, $E_1$ and $E_2$. If $d_{122} \neq 0$ and $d_{112} \neq 0$, then all three equilibriums are hyperbolic and the conclusions regarding their stability are given in Table~\ref{Table5}. 
	
\begin{table}[h!]
	\centering
		\begin{tabular}{||>{\centering\arraybackslash}m{0.3cm}||>{\centering\arraybackslash}m{2.1cm}|>{\centering\arraybackslash}m{2.2cm}|>{\centering\arraybackslash}m{2.7cm}|>{\centering\arraybackslash}m{2.1cm}||}
		\hline
		  eq. & $d_{12}=0$ \newline $d_{112}>0$ \newline $d_{122}>0$  & $d_{12}=0$ \newline $d_{112}\, d_{122}<0$ & $d_{12}=0$ \newline $d_{112}\, d_{122}=0$ \newline ${d_{112}}^2+{d_{122}}^2 \neq 0$ & $d_{12}=0$ \newline  $d_{112}<0$ \newline $d_{122}<0$\\
		\hline\hline
		  $E_0$ & U (u.\ node) & \multirow{3}{1.9cm}{\centering not possible cases} & \multirow{3}{1.9cm}{\centering not possible cases} & U (u.\ node)\\	\cline{1-2} 	 \cline{5-5}
		  $E_1$ & U (saddle) &   & &  AS (s.\ node)    \\	\cline{1-2} 	 \cline{5-5}	
		  $E_2$ & AS (s.\ node) &  & &  U (saddle)   \\ 
		\hline
		\end{tabular}
		\caption{Stability analysis of the equilibriums of the system (\ref{2d}).}
		\label{Table5}
	\end{table}

\noindent Here the cases for $d_{12}=0$ and $d_{112}\, d_{122}<0$ are impossible and the cases for $d_{12}=0$, $d_{112}\, d_{122}=0$ and ${d_{112}}^2+{d_{122}}^2 \neq 0$ are also impossible. The proofs are similar to those in Section~\ref{nzero2} and so we omit them here.

\noindent We note that this section does not provide us with any qualitatively new phase portraits of the dynamical system (\ref{2d}), since we already had these qualitatively equal phase portraits before. For $d_{12}=0$, $d_{112}>0$, $d_{122}>0$, we have the same phase portrait as for $d_{12}>0$, $d_{112}>0$, $d_{122}>0$ and $d_{12}<0$, $d_{112}>0$, $d_{122}>0$ from Section~\ref{nzero12}. Furthermore, for $d_{12}=0$, $d_{112}<0$, $d_{122}<0$, we have the same phase portrait as for the $d_{12}>0$, $d_{112}<0$, $d_{122}<0$ and $d_{12}<0$, $d_{112}<0$, $d_{122}<0$ from Section~\ref{nzero12}.

\subsubsection{Case \boldmath{$d_{122} = 0$} and \boldmath{$d_{112} = 0$}}\label{zero12}

If both minors are equal to zero, the system (\ref{2ds}) has infinitely many solutions $\left(x_1,x_2\right)=\left(\left(b_1-a_{12} \alpha \right)/a_{11},\alpha \right)$, where $\alpha \geq 0$ and $b_1 \geq a_{12} \alpha$. Consequently, the system (\ref{2d}) has infinitely many equilibriums, $E_0$ and $E^{\alpha}=\left(\left(b_1-a_{12} \alpha \right)/a_{11},\alpha \right)$, where $0 \leq \alpha \leq b_1/a_{12}$. We note that the equilibriums $E^{\alpha}$ are non-isolated and that for $\alpha=0$ we conclude $E^0=E_1$, while for $\alpha=b_1/a_{12}$ we obtain $E^{b_1/a_{12}}=E^{b_2/a_{22}}=E_2$. Therefore we exclude these values for $\alpha$ in the following discussion.

\noindent The Jacobian matrix calculated for $E^{\alpha}$ is given by
	\begin{equation}\label{jacalpha}
		\begin{array}{ll}
		\smallskip
		\displaystyle J\left(E^{\alpha}\right)&=\left[\begin{matrix}
											-\left(b_1-a_{12}\alpha\right) & -\dfrac{\left(b_1-a_{12}\alpha\right)a_{12}}{a_{11}}\\
											-a_{21}\alpha & \dfrac{a_{11}\left(b_2-2a_{22}\alpha\right)-a_{21}\left(b_1-a_{12}\alpha\right)}{a_{11}}
											\end{matrix}\right].					\end{array}	
	\end{equation}	
From (\ref{jacalpha}), using $d_{112}=0$ and $d_{12}=0$, we deduce that
	\begin{equation}\label{tealpha}
		\begin{array}{l}
		\displaystyle \lambda_1 \lambda_2 = \det\left( J\left(E^{\alpha}\right)\right)=0,\\
		
		\smallskip
		\displaystyle \lambda_1+\lambda_2 = \tr \left( J\left(E^{\alpha}\right)\right)=a_{12}\alpha-a_{22}\alpha-b_1,		
		\end{array}	
	\end{equation}
from which we conclude that $\displaystyle \lambda_2=-\left(b_1-a_{12}\alpha\right)-a_{22}\alpha$, $\lambda_2<0$, $\displaystyle \lambda_1=0$, since $\left(b_1-a_{12}\alpha\right)/a_{11}=x_1$, $x_1>0$ and $\alpha=x_2$, $x_2>0$. Thus, $E^{\alpha}$ is nonhyperbolic equilibrium for every $0<\alpha<b_1/a_{12}$.
$E^{\alpha}$ is a nonhyperbolic equilibrium for every $0 \leq \alpha \leq b_1/a_{12}$, i.e.\ $x_2=\left(b_1-a_{11}x_1\right)/a_{12}$ is a line of equilibriums. $E_0$ is a hyperbolic equilibrium, which is an unstable node, just as in the previous sections. 
The phase portrait is shown in Figure~\ref{f5}.

\begin{figure}[h]
	\centering 
	\captionsetup{justification=centering}
	\includegraphics[width=0.44\textwidth]{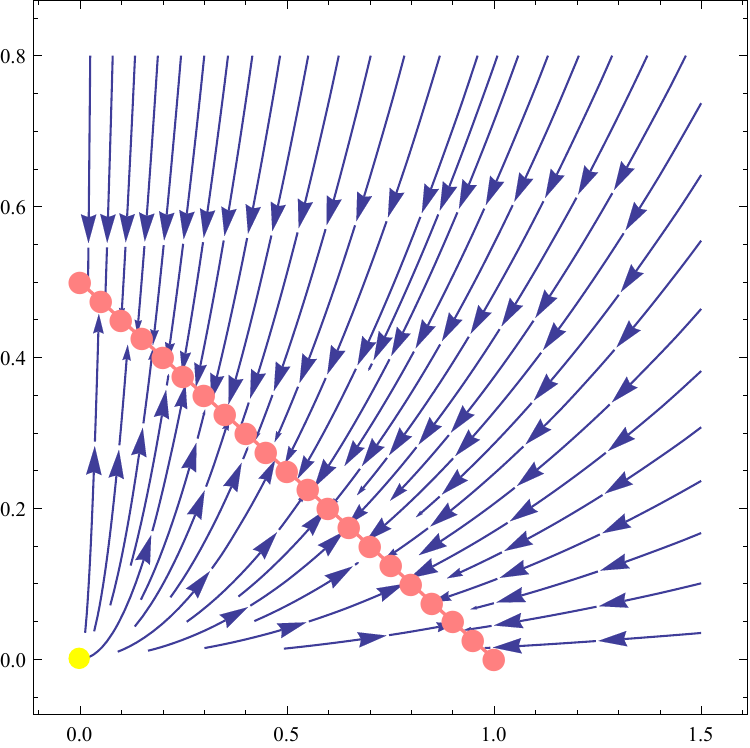}	 
	\caption{Case $d_{12}=d_{112}=d_{122}=0$. \\$\dot{x_1}=x_1 \left(1-x_1-2x_2\right)$, $\dot{x_2}=x_2 \left(2-2x_1-4x_2\right)$.}
		 \label{f5}
	\end{figure} 
	
\noindent The equilibrium $E_0$ at the origin is unstable, while the equilibriums $E^{\alpha}$, where $0 \leq \alpha \leq b_1/a_{12}$, are non-isolated.

\section{Necessary and sufficient conditions for stability of equilibriums}\label{theorems}

Section 2 shows that there are nine qualitatively different phase portraits of
the dynamical system (1.1) which are shown in Figures~\ref{f1}-\ref{f5}, when we consider the stability of equilibriums of the system (1.1) on their full neighbourhoods. The signs of $d_{12}$, $d_{112}$ and $d_{122}$ as well as the equilibriums of the dynamical system (\ref{2d}) are shown in Table~\ref{Table6}. 
	\begin{table}[h!]
	\centering
		\begin{tabular}{||>{\centering\arraybackslash}m{0.8cm}||>{\centering\arraybackslash}m{1.3cm}|>{\centering\arraybackslash}m{1.4cm}|>{\centering\arraybackslash}m{1.4cm}||>{\centering\arraybackslash}m{0.5cm}||>{\centering\arraybackslash}m{0.5cm}||>{\centering\arraybackslash}m{0.5cm}||>{\centering\arraybackslash}m{0.5cm}||>{\centering\arraybackslash}m{0.5cm}||}
		\hline
		Serial No. & $d_{12}$ & $d_{112}$ & $d_{122}$ & $E_0$ & $E_1$ & $E_2$ & $E_{12}$ & $E_{\alpha}$\\
		\hline
		\hline
		1 & $d_{12}>0$ & $d_{112}>0$ & $ d_{122}<0$ & U & U & U & AS & / \\
		\hline
		2 & $d_{12}>0$ & $d_{112}>0$ & $ d_{122}=0$ & U & U & SS & / & / \\
		\hline
		3 & $d_{12}>0$ & $d_{112}>0$ & $ d_{122}>0$ & U & U & AS & / & / \\
		\hline
		4 & $d_{12}>0$ & $d_{112}=0$ & $ d_{122}<0$ & U & SS & U & / & / \\
		\hline
		5 & $d_{12}>0$ & $d_{112}<0$ & $ d_{122}<0$ & U & AS & U & / & / \\
		\hline
		6 & $d_{12}<0$ & $d_{112}=0$ & $ d_{122}>0$ & U & U & AS & / & / \\
		\hline
		7 & $d_{12}<0$ & $d_{112}<0$ & $ d_{122}=0$ & U & AS & U & / & / \\
		\hline
		8 & $d_{12}<0$ & $d_{112}<0$ & $ d_{122}>0$ & U & AS & AS & U & / \\
		\hline
		9 & $d_{12}=0$ & $d_{112}=0$ & $ d_{122}=0$ & U & NI & NI & / & NI \\
		\hline
		\end{tabular}
		\caption{Stability analysis of equilibriums of the dynamical system (\ref{2d}).}
		\label{Table6}
	\end{table}
After analyzing Figures~\ref{f1}-\ref{f5} and the data in Table~\ref{Table6}, we conclude that the phase portraits under serial numbers 2, 3 and 6 that are shown in Figures~\ref{f3a}, \ref{f1b} and \ref{f4b} have qualitatively the same dynamical properties in $\left[0, +\infty\right)^2$, i.e.\ that in $\left[0, +\infty\right)^2$ the equilibrium $E_1$ is an unstable equilibrium and that the equilibrium $E_2$ is asymptotically stable.

\noindent Moreover, the phase portraits under serial numbers 4, 5 and 7 that are shown in Figures~\ref{f4a}, \ref{f2a} and \ref{f3b} have qualitatively the same dynamical properties in $\left[0, +\infty \right)^2$, i.e.\ the equilibrium $E_1$ is asymptotically stable and the equilibrium $E_2$ is unstable. 

\noindent In summary, we have nine qualitatively different phase portraits of the dynamical system (\ref{2d}) when we consider the stability of equilibriums of the system (\ref{2d}) on their full neighbourhoods, while we have five qualitatively different phase portraits when we analyze the stability of its equilibriums on $\left[0, +\infty \right)^2$. We note that it is interesting from mathematical point of view  to analyze the stability of the equilibriums of the system (\ref{2d}) on their full neighbourhoods, although it lacks an appropriate physical interpretation since the density of the species cannot be negative. 

\noindent To summarize all the results from the Section~\ref{sectionstability}, we formulate the following theorems. 

\begin{theorem}\label{thrE1SN}
The equilibrium $E_1\left(b_1/a_{11},0\right)$ of the system (\ref{2d}) is asymptotically stable on $\left[0, +\infty \right)^2$ if and only if one of the following conditions is satisfied:
\begin{itemize}
\item[(1)]  $d_{112}<0$ and $d_{122}<0$ (regardless of the $\mathrm{sgn}\left(d_{12}\right)$);
\item[(2)] $d_{12}<0$ and $d_{112}<0$ and $d_{122}=0$;
\item[(3)] $d_{12}>0$ and $d_{112}=0$ and $d_{122}<0$;
\item[(4)] $d_{12}<0$ and $d_{112}<0$ and $d_{122}>0$.
\end{itemize}
\end{theorem}
Theorem~\ref{thrE1SN} is an extension of Theorem 2.1. from \cite{zeeman1995} for $n=2$, where it states (with a different notation than the one used here) that if $d_{122}<0$ and $d_{112}<0$, then $E_1\left(b_1/a_{11},0\right)$ is asymptotically stable on $\left(0, +\infty \right)^2$. 

\noindent Now we formulate a similar theorem for the equilibrium $E_2$.
\begin{theorem}\label{thrE2SN}
The equilibrium $E_2\left(0, b_2/a_{22}\right)$ of the system (\ref{2d}) is asymptotically stable on $\left[0, +\infty \right)^2$ if and only if one of the following conditions is satisfied:
\begin{itemize}
\item[(1)]  $d_{112}>0$ and $d_{122}>0$ (regardless of the $\mathrm{sgn}\left(d_{12}\right)$);
\item[(2)] $d_{12}>0$ and $d_{112}>0$ and $d_{122}=0$;
\item[(3)] $d_{12}<0$ and $d_{112}=0$ and $d_{122}>0$;
\item[(4)] $d_{12}<0$ and $d_{112}<0$ and $d_{122}>0$.
\end{itemize}
\end{theorem}

\noindent As a consequence of Theorem~\ref{thrE1SN} and Theorem~\ref{thrE2SN} and the fact that $E_0$ is an unstable node (unstable equilibrium) regardless of the signs of $d_{12}$, $d_{112}$ and $d_{122}$, we formulate the following two theorems.

\begin{theorem}\label{corE1S}
The equilibrium $E_1\left(b_1/a_{11},0\right)$ of the system (\ref{2d}) is unstable if and only if one of the conditions $\left(1\right)$, $\left(2\right)$, $\left(3\right)$ from Theorem~\ref{thrE2SN} is satisfied, or $d_{12}>0$ and $d_{112}>0$ and $d_{122}<0$ is satisfied.
\end{theorem}

\begin{theorem}\label{corE2S}
The equilibrium $E_2\left(0, b_2/a_{22}\right)$ of the system (\ref{2d}) is unstable if and only if one of the conditions $\left(1\right)$, $\left(2\right)$, $\left(3\right)$ from Theorem~\ref{thrE1SN} is satisfied, or $d_{12}>0$ and $d_{112}>0$ and $d_{122}<0$ is satisfied.
\end{theorem}

\noindent Moreover, the following can be concluded. 

\begin{theorem}\label{nofp}
The system (\ref{2d}) has no equilibriums in $\left(0, +\infty \right)^2$ if and only if one of the conditions $\left(1\right)$, $\left(2\right)$, $\left(3\right)$ from Theorem~\ref{thrE1SN} or one of the conditions $\left(1\right)$, $\left(2\right)$, $\left(3\right)$ from Theorem~\ref{thrE2SN} is satisfied.
\end{theorem}
Theorem~\ref{nofp} is an extension of Lemma 7.2. from \cite{zeeman1995} for $n=2$, which states (with a different notation than the one used here) that if the system (\ref{2d}) satisfies the inequalities $d_{122}<0$ and $d_{112}<0$, then there is no equilibrium in $\left(0, +\infty \right)^2$. The next theorem considers the stability of equilibrium $E_{12}$.
\begin{theorem}\label{thrE12}
The equilibrium $E_{12}\left(-d_{122}/d_{12},d_{112}/d_{12}\right)$ of the system (\ref{2d}) is:
\begin{itemize}
\item[(1)] asymptotically stable on $\left(0, +\infty \right) ^2$ if and only if $d_{12}>0$ and $d_{112}>0$ and $d_{122}<0$;
\item[(2)] unstable if and only if $d_{12}<0$ and $d_{112}<0$ and $d_{122}>0$.
\end{itemize}
\end{theorem}

\noindent Proofs for all theorems can be found in discussions in the Section~\ref{sectionstability} and so are omitted here. We note that in \cite{florian} are considered only the cases when appropriate determinants are non-zero, while here we investigated those cases too.

\section{Bifurcation analysis}\label{sectionbifurcation}

As a consequence of Section~\ref{sectionstability}, there are five different phase portraits of the system (\ref{2d}) in $\left[0, +\infty \right)^2$, which can be found in Table~\ref{Table6} under serial numbers 1, 2, 4, 8 and 9. They are already shown in Figures~\ref{f1a}, \ref{f3a}, \ref{f4a}, \ref{f2b} and \ref{f5} respectively. We note that it is interesting from mathematical point of view to analyze the stability of the equilibriums of the system (\ref{2d}) on $\mathrm{R}^2$, although it lacks an appropriate physical interpretation. In that case, four transcritical bifurcations can occur. 

\noindent Comparing cases 1, 2 and 3 from Table~\ref{Table6} but with the stability of equilibriums analyzed on $\mathrm{R}^2$, we see that for $d_{122}<0$ (case 1) all equilibriums are hyperbolic and that $E_0$, $E_1$ and $E_2$ are an unstable equilibriums, while $E_{12}$ is an asymptotically stable equilibrium. For $d_{122}=0$ (case 2), nonhyperbolic equilibrium $E_2$ collides with $E_{12}$ and becomes a semi-stable equilibrium, while the stability of other two equilibriums remains the same. Furthermore, for $d_{122}>0$ (case 3), the equilibrium $E_{12}$ appears again, but this time as an unstable equilibrium in the second quadrant, while the equilibrium $E_2$ is now an asymptotically stable equilibrium, i.e.\ the equilibrium $E_2$ have swapped stability with the equilibrium $E_{12}$ at $d_{122}=0$, while the stability of other equilibriums have remained the same.

\noindent We notice similar situation when we compare cases 1, 4 and 5 from Table~\ref{Table6}, again analyzing the system (\ref{2d}) on $\mathrm{R}^2$. More precisely, the equilibrium $E_1$ have swapped stability with $E_{12}$ (which transitioned from the first quadrant to the fourth quadrant) at $d_{112}=0$. Furthermore, when we compare cases 8 and 6 from Table~\ref{Table6} and the case $d_{12}<0$, $d_{112}>0$ and $d_{122}>0$ from Table~\ref{Table1} (again on $\mathrm{R}^2$), we easily obtain that the equilibrium $E_1$ have swapped stability with the equilibrium $E_{12}$ (which transitioned from the first quadrant to the fourth quadrant) at $d_{112}=0$ if $a_{11}x_1 + a_{22}x_2>0$ $\left(a_{22}d_{112}<a_{11}d_{122}\right)$. Finally, when we compare cases 8 and 7 from Table~\ref{Table6} and the case $d_{12}<0$, $d_{112}<0$ and $d_{122}<0$ from Table~\ref{Table1}, we derive that the equilibrium $E_2$ have swapped stability with the equilibrium $E_{12}$ (which transitioned from the first quadrant to the second quadrant) at $d_{122}=0$, again if $a_{11}x_1 + a_{22}x_2>0$.

\section{Conclusion}

We have analyzed the stability of equilibriums and bifurcations of the dynamical system (\ref{2d}) for various signs of the main and minor determinants of this system, $d_{12}$, $d_{112}$ and $d_{122}$. Of the total of twenty-seven different combinations of the signs of the determinants $d_{12}$, $d_{112}$ and $d_{122}$, thirteen were meaningful, while nine of them were different from each other on full neighbourhood of the equilibriums and their phase portraits are shown in Figures~\ref{f1}-\ref{f5}. Only five of these nine cases had different phase portraits on $\left[0, +\infty\right)$, which are shown in Figures~\ref{f1a}, \ref{f3a}, \ref{f4a}, \ref{f2b} and \ref{f5}. Four of these five phase portraits can also be found in \cite{florian}, while the last one is new and very specific as it has infinite equilibriums. All results are presented in Tables~\ref{Table1}-\ref{Table6}. In addition, we have improved several results from \cite{zeeman1995}, mainly in terms of necessary and sufficient conditions for the equilibriums $E_1$, $E_2$ and $E_{12}$ to be asymptotically stable or unstable. We noticed four transcritical bifurcations among thirteen meaningful and different phase portraits on $\mathrm{R}^2$.

	School of Electrical Engineering, University of Belgrade,
	Bulevar kralja Aleksandra 73, 11120 Belgrade, Serbia\\
	danijela@etf.bg.ac.rs\\
	\\
	Faculty of Mathematics, University of Belgrade,
	Studentski trg 16, 11158 Belgrade, Serbia \\
	marija.mikic@matf.bg.ac.rs\\
\end{document}